\pdfoutput=1
\RequirePackage{ifpdf}
\ifpdf %We are running pdfTeX in pdf mode
\documentclass[pdftex]{sigma}
\else
\documentclass{sigma}
\fi

\begin{document}

\allowdisplaybreaks

\renewcommand{\PaperNumber}{023}

\FirstPageHeading

\ShortArticleName{Non-Schlesinger Isomonodromic Deformations of Fuchsian Systems and Middle Convolution}

\ArticleName{Non-Schlesinger Isomonodromic Deformations\\
of Fuchsian Systems and Middle Convolution}

\Author{Yulia BIBILO~$^\dag$ and Galina FILIPUK~$^\ddag$}

\AuthorNameForHeading{Yu.~Bibilo and G.~Filipuk}

\Address{$^\dag$~Department of Theory of Information Transmission and Control, Institute for Information\\
\hphantom{$^\dag$}~Transmission Problems, Russian Academy of Sciences, Bolshoy Karetny per.~19,\\
\hphantom{$^\dag$}~Moscow, 127994, Russia}
\EmailD{\href{mailto:y.bibilo@gmail.com}{y.bibilo@gmail.com}}
\URLaddressD{\url{http://www.iitp.ru/en/users/2131.htm}}

\Address{$^\ddag$~Faculty of Mathematics, Informatics and Mechanics, University of Warsaw,\\
\hphantom{$^\dag$}~Banacha~2, Warsaw, 02-097, Poland}
\EmailD{\href{mailto:filipuk@mimuw.edu.pl}{filipuk@mimuw.edu.pl}}
\URLaddressD{\url{http://www.mimuw.edu.pl/~filipuk/}}

\ArticleDates{Received November 20, 2014, in f\/inal form March 04, 2015; Published online March 13, 2015}

\Abstract{The paper is devoted to non-Schlesinger isomonodromic deformations for resonant Fuchsian systems.
There are very few explicit examples of such deformations in the literature.
In this paper we construct a~new example of the non-Schlesinger isomonodromic deformation for a~resonant Fuchsian system of
order~5 by using middle convolution for a~resonant Fuchsian system of order~2.
Moreover, it is known that middle convolution is an operation that preserves Schlesinger's deformation equations for
non-resonant Fuchsian systems.
In this paper we show that Bolibruch's non-Schlesinger deformations of resonant Fuchsian systems are, in general, not
preserved by middle convolution.}


\Classification{34M56; 44A15}

\section{Introduction}

Dettweiler and Reiter's algebraic analogue~\cite{D, DR, DR2} of Katz' middle convolution~\cite{Katz} is a~certain
transformation of Fuchsian systems which preserves an index of rigidity.
Middle convolution is related to the Euler transformation of solutions of the Fuchsian systems.
There have been numerous studies on middle convolution in recent years, including applications to special functions
(e.g.,~\cite{DR_Painleve, GF_Painleve, GF_Heun, GF_hyper, Reiter, Takemura}), extensions
to irregular systems
(e.g.,~\cite{Hiroe, Kawakami, Takemura_irreg}) and various others (see, for instance,~\cite{Arinkin, DeligneSimson,Silva,D_motives, D,
D_Sabbah, Haraoka, Oshima, Oshima_frac, P, Simpson, Yamakawa, Yamakawa2}).

The theory of isomonodromic deformations for Fuchsian systems was developed~by
L.~Schle\-singer~\cite{Schlesinger,Schlesinger2}, L.~Fuchs~\cite{F}, R.~Garnier~\cite{G,G2} (see
also~\cite{GaussPainleve, Jimbo, Malgrange} and others).
Isomonodromy deformations mean that a~monodromy group of the Fuchsian system does not depend on the parameters (the
location of poles).
It is known that in this case the residue matrices satisfy the so-called Schlesinger equation (with certain assumptions
to be discussed below).
Isomonodromic deformations are, in particular, important in the theory of special functions (especially the Painlev\'e
transcendents), random matrices, the Riemann--Hilbert problems.
They are closely connected to the integrability and the Painlev\'e property.

\looseness=-1
It is known that there exist non-Schlesinger deformations~\cite{Bol7, Bol3, Bol1}.
The residue matrices of such deformations do not satisfy the Schlesinger equation, but for any values of the pa\-ra\-me\-ters
the family of Fuchsian systems has the same monodromy representation.
The theory of non-Schlesinger isomonodromic deformations completes the theory of the isomonodromic deformations.

The study of the resonant Fuchsian systems appeared to be very fruitful in the areas related to the analytic theory of
linear dif\/ferential equations and led to some unexpected results (see~\cite{Il_Yak}).
In particular, non-Schlesinger isomonodromic deformations appeared in the following problems.
A.A.~Bolibruch considered a~general isomonodromic dif\/ferential 1-form in the context of isomonodromic conf\/luence of
Fuchsian singularities.
The question of isomonodromic conf\/luence was posed by V.I.~Arnold
and it was included in his list of problems (see \cite[Pr.~1984-7, Pr.~1987-12]{Arnold}).
In particular, one of the problems was to get a~system with irregular singularities as a~limit of a~Fuchsian
isomonodromic family when its singularities coalesce.
A.A.~Bolibruch proved that it is impossible using the theory of non-Schlesinger isomonodromic deformations.
The general isomonodromic dif\/ferential 1-form was crucial in the solution of this problem.
In addition, the existence of examples of non-Schlesinger deformations essentially answers the question posed~by
Y.~Sibuya in~\cite{Sibuya}.
He asked whether it is necessary for the connection matrices of the isomo\-nodromic family with the Fuchsian and irregular
singularities to be constant.
The answer is that it it not necessary (see Def\/inition~\ref{def_isom} below).

There exist some further studies and extensions of the notion of isomonodromic deformations of   the Fuchsian systems
(see, for instance,~\cite{MitschiSinger, Pober}).
In this paper we deal with the def\/inition of isomonodromic deformations as in Def\/inition~\ref{def_isom},
Theorems~\ref{th_isom} and~\ref{th 2.2}.

In~\cite{HaraokaFilipuk} it was shown that the deformation equations for Fuchsian equations are preserved by middle
convolution.
In particular, the so-called Hitchin systems~\cite{Hitchin} obtained from the Schlesinger systems are invariant under
middle convolution.
One of the aims of this paper is to study non-Schlesinger deformations and their behaviour under middle convolution.

Usually it is very dif\/f\/icult to f\/ind explicit examples of non-Schlesinger deformations.
The f\/irst example of the non-Schlesinger isomonodromic deformation which cannot be reduced to the Schlesinger one was
constructed by A.A.~Bolibruch in~\cite{Bol8}.
To date, very few examples exist~\cite{Bibilo_Thesis,KatVolok,Kitaev}.
For instance, to our knowledge, there are no explicit examples of non-Schlesinger deformations for resonant Fuchsian
systems with higher order poles in the ``non-Schlesinger'' part of the isomonodromic dif\/ferential form, which is def\/ined
below.
All of the known examples are for Fuchsian systems of orders~2 or~3.
Therefore, a~construction of explicit examples for systems of higher order is interesting and new in its own right.

The paper is organized as follows.
In the f\/irst section we recall the theory of Schlesinger and non-Schlesinger deformations following~\cite{Bol3,Bol1,
GaussPainleve}.
Next we brief\/ly summarize the algebraic construction of middle convolution following~\cite{D, DR, DR2}.
We present a~new explicit example of the non-Schlesinger isomonodromic deformation of a~Fuchsian system of order~5.
Finally we present an explicit example which shows that non-Schlesinger isomonodromic deformations are not preserved~by
middle convolution in general.
We also discuss a~number of open questions and problems.
Some of them are raised by the examples in the current paper.

\section{Isomonodromic deformations}

Let us consider a~system of~$p$ linear dif\/ferential equations on the Riemann sphere
\begin{gather}
\label{syst0}
\frac{dy}{dz}=\left(\sum\limits_{i=1}^n\frac{A_i^0}{z-a_i^0}\right)y,
\qquad
\sum\limits_{i=1}^n A_i^0=-A^0_{n+1}
\end{gather}
with distinct singular points $a^0_1, \ldots, a^0_n \in \mathbb{C}$, $a^0_{n+1}=\infty$.
Here  $y(z)\in\mathbb{C}^p$ is unknown.
System~\eqref{syst0} is a~Fuchsian system and its important characteristic is a~monodromy group.
The monodromy representation is def\/ined as follows.
Take an arbitrary non-singular point $z_0\in \mathbb{C}$.
Let $Y(z)$ be a~fundamental matrix solution of system~\eqref{syst0} and $\tilde{Y}(z)$ be its analytic continuation
along a~loop~$\gamma$ which starts and ends at $z_0$.
These fundamental matrices are related by $Y(z)=\tilde{Y}(z)G_{\gamma}$, where the matrix $G_{\gamma}$ is constant with
respect to~$z$ and depends only on the homotopy class $[\gamma]$ of the loop~$\gamma$.
The mapping $\gamma\to G_{\gamma}$ def\/ines a~representation of the fundamental group into the group of non-degenerate
matrices, i.e.,
\begin{gather}
\label{monodr0}
\chi^0: \ \pi_1 \big(T^0,z_0\big) \rightarrow \mathbb{GL}(p, \mathbb{C}),
\qquad
T^0=\bar{\mathbb{C}}\backslash \cup_{i=1}^{n+1}\big\{a^0_i\big\}.
\end{gather}
The generators $\gamma_1,\ldots,\gamma_n,\gamma_{n+1}$ of the fundamental group (the loops around singularities
$\{a^0_1$, $\ldots, a^0_n,a^0_{n+1}\}$) satisfy $\gamma_1 \cdots \gamma_n \gamma_{n+1}=e$.

Now consider a~family of Fuchsian systems
\begin{gather}
\label{syst}
\frac{dy}{dz}=\left(\sum\limits_{i=1}^n\frac{A_i(a)}{z-a_i}\right)y,
\qquad
\sum\limits_{i=1}^n A_i(a)=-A_{n+1}(a),
\end{gather}
where the location of singularities varies, which is denoted by a~parameter~$a$, i.e., $a=(a_1,\ldots,a_n)$.
Moreover, we assume that the initial system~\eqref{syst0} is included into the isomonodromic family of~\eqref{syst} with
$A_i(a^0)=A_i^0$, $a^0=(a_1^0, \ldots, a_n^0)$, $i=1,\ldots,n$.
Similarly to~\eqref{monodr0} we can def\/ine the monodromy representation of system~\eqref{syst} for any $a\in
D(a^0)\backslash \cup_{i,j=1, i\neq j}^n \{a_i=a_j\}$~by
\begin{gather*}
\chi_a: \ \pi_1 (T_a, z_0) \rightarrow \mathbb{GL}(p, \mathbb{C}),
\qquad
T_a=\bar{\mathbb{C}}\backslash \cup_{i=1}^{n+1} \{a_i\}.
\end{gather*}
Here $D(a^0)$ is a~small open disk centered at $a^0$.

\begin{definition}[\cite{Bol2,Bol3, Bol1}]\label{def_isom}
The family of Fuchsian systems \eqref{syst} is called isomonodromic if the monodromy representation $\chi_a$ coincides
with the monodromy representation $\chi^0$ of the initial system \eqref{syst0} for any $a\in D(a^0)\backslash
\cup_{i, j=1, i\neq j}^n \{a_i=a_j\}$.
\end{definition}
Def\/inition~\ref{def_isom} means that
there exists a~fundamental solution $Y(z,a)$ (called an isomonodromic matrix) def\/ining the
same monodromy matrices for any $a\in D(a^0)\backslash \cup_{i, j=1, i\neq j}^n \{a_i=a_j\}$.
The isomonodromic family \eqref{syst} is also called the isomonodromic deformation.
There exists a~ca\-no\-nical isomorphism between the fundamental groups $\pi_1(T^0,z_0)$ and $\pi_1(T_a,z_0)$ (see, for
instance,~\cite{Bol3}).
Therefore, the equality between the monodromy representations $\chi_a$ and $\chi^0$ is def\/ined correctly.

The following theorem was formulated and proved for $(2\times 2)$-families of Fuchsian systems in~\cite{FN}.
A.A.~Bolibruch generalized it for families of any order.

\begin{theorem}[\cite{Bol2,Bol3, Bol1}]\label{th_isom}
The family of Fuchsian systems \eqref{syst} is isomonodromic if and only if there exists a~matrix-valued differential
$1$-form~$\omega$ on $\mathbb{C}\times D(a^0)\backslash \cup_{i=1}^n\{z-a_i=0\}$ such that
\begin{enumerate}\itemsep=0pt
\item[$i)$] $\omega=\sum\limits_{i=1}^n \frac{A_i(a)}{z-a_i}dz$ for any fixed $a\in D(a^0)$;

\item[$ii)$]
$d\omega=\omega \wedge \omega $.
\end{enumerate}
\end{theorem}

Def\/inition~\ref{def_isom} was proposed by A.A.~Bolibruch (see, for instance,~\cite{Bol1}).
It extends the notion of the isomonodromic deformation which was used before.
The Schlesinger deformations are the most known and well-studied isomonodromic deformations in the literature.
They are given by the dif\/ferential 1-form
\begin{gather}
\label{schl}
\omega_{\rm Schl}=\sum\limits_{i=1}^n\frac{A_i(a)}{z-a_i}d(z-a_i).
\end{gather}
One can verify directly that condition ii) of Theorem~\ref{th_isom} is equivalent to
\begin{gather*}
dA_i(a)=-\sum\limits_{j=1, j\neq i}^n \frac{[A_i(a), A_j(a)]}{a_i-a_j}d(a_i-a_j),
\end{gather*}
which is called the Schlesinger equation.
It is easily seen that the residue matrices satisfy the so-called Schlesinger equations
\begin{gather}
\label{Schl eqs}
\frac{\partial A_i(a)}{\partial a_i}=-\sum\limits_{k\neq i}\frac{[A_i(a),A_k(a)]}{a_i-a_k},
\qquad
\frac{\partial A_j(a)}{\partial a_i}=\frac{[A_i(a),A_j(a)]}{a_i-a_j},
\quad
j\neq i.
\end{gather}
Note that in~\cite{Hitchin} N.~Hitchin derives a~new system of equations from the Schlesinger equations~\eqref{Schl eqs}
for traces of products of matrices and their commutators ${\rm tr} ([A_i,A_j]A_k)$.

The fundamental matrix $Y_{\rm Schl}(z, a)$ of the isomonodromic family def\/ined by the dif\/ferential 1-form
\eqref{schl} satisf\/ies
\begin{gather*}
Y_{\rm Schl}(\infty, a) \equiv C,
\end{gather*}
where~$C$ is a~constant non-degenerated matrix. Indeed, the 1-form of the deformation is def\/ined (see~\cite{Bol1}) by
using the fundamental matrix $Y(z)$ of the family as follows:
\begin{gather*}
\omega=dY\cdot Y^{-1},
\end{gather*}
and, hence,
\begin{gather*}
d_a Y_{\rm Schl}(\infty, a)Y_{\rm Schl}^{-1}(\infty, a)=-\sum\limits_{i=1}^n \frac{A_i(a)}{z-a_i}d(a_i)|_{z=\infty}\equiv 0.
\end{gather*}

\begin{definition}
%\label{def_norm}
A~deformation of the Fuchsian system is called normalized if the choice of the isomonodromic fundamental matrix of the
system is f\/ixed by $Y(\infty, a) \equiv C$, where the matrix~$C$ is constant and non-degenerate.
\end{definition}

Let us return to the def\/inition of the monodromy representation of the Fuchsian system \eqref{syst0}.
Take another fundamental matrix $X(z)$ of system~\eqref{syst0}.
Two fundamental matrices~$Y(z)$ and~$X(z)$ are related by $X(z)=Y(z)M$, where~$M$ is a~constant non-degenerate matrix.
The analytic continuation~$\tilde{X}(z)$ of the fundamental matrix $X(z)$ along a~loop~$\gamma$  satisf\/ies
$\tilde{X}(z)=X(z)\tilde{G}_{\gamma}$.
Then the monodromy matrices~$G_{\gamma}$ and~$\tilde{G}_{\gamma}$ are related~by
\begin{gather*}
\tilde{G}_{\gamma}=M^{-1} G_{\gamma} M.
\end{gather*}
Therefore, the Fuchsian system \eqref{syst0} def\/ines not only a~unique monodromy representation $\chi^0$ but
a~conjugacy (by a~constant matrix) class of representations.

To summarise,
Def\/inition~\ref{def_isom} proposed by A.A.~Bolibruch is a~natural and well-grounded extension of the notion of the normalized
isomonodromic deformation since there is no any requirement that the monodromy representation for any value of the
parameter is def\/ined with respect to the fundamental matrix f\/ixed by some principle (see~\cite{Anosov1} for more details).

An arbitrary isomonodromic deformation is not necessarily the Schlesinger deformation~\eqref{schl}.
Let us consider a~family of Fuchsian systems with the fundamental matrix
\begin{gather*}
Y(z,a)=\Gamma(a)Y_{\rm Schl}(z,a),
\end{gather*}
where~$\Gamma(a)$ is a~holomorphically invertible matrix, $\Gamma(a_0)=I$.
In this case the dif\/ferential 1-form $\omega=dY(z,a)Y^{-1}(z,a)$ is given~by
\begin{gather}
\label{simple_nonschl}
\omega=\sum\limits_{i=1}^n\frac{A'_i(a)}{z-a_i}d(z-a_i)+\sum\limits_{k=1}^n \gamma_k(a)da_k,
\end{gather}
where the coef\/f\/icients are given by~\cite{Bol3, Bol1}
\begin{gather*}
A'_i(a)=\Gamma(a)A_i(a)\Gamma^{-1}(a),
\qquad
\gamma_k(a)=\frac{\partial \Gamma(a)}{\partial a_k} \Gamma^{-1}(a).
\end{gather*}
This isomonodromic deformation is not normalized.
It is clear that this deformation is reduced to the Schlesinger deformation~by
\begin{gather*}
Y_{\rm Schl}(z,a)=\Gamma^{-1}(a)Y(z,a).
\end{gather*}
Note that one can include initial system \eqref{syst0} into two isomonodromic deformations def\/ined by~$\omega_{\rm Schl}$
and~$\omega$.

However, there exist isomonodromic deformations given by dif\/ferential 1-forms dif\/ferent from~\eqref{schl} and~\eqref{simple_nonschl}.
In~\cite{Bol1} A.A.~Bolibruch gave examples of such deformations and obtained a~general form of the isomonodromic
deformation.

\begin{definition}[\cite{Bol2, Bol1}]
Let $\lambda^i_1, \ldots, \lambda^i_p$ be the eigenvalues of the matrix~$A_i$ of the Fuchsian system~\eqref{syst}.
A~singular point~$a_i$ is called resonant if there exist at least two distinct eigenvalues of~$A_i$ such that their
dif\/ference is a~natural number.
The number
\begin{gather*}
r_i= \max_{k \neq j,|\lambda^i_k-\lambda^i_j|\in \mathbb{N}}\big|\lambda^i_k-\lambda^i_j\big|
\end{gather*}
is called the maximal~$i$-resonance of the system.
\end{definition}

\begin{theorem}[\cite{Bol2, Bol1}]\label{th 2.2}
Any matrix-valued differential $1$-form~$\omega$ on $\bar{\mathbb{C}}\times D(a^0)\backslash \cup_{i=1}^n\{z-a_i=0\}$
which defines the isomonodromic deformation of \eqref{syst}  is given~by
\begin{gather*}
\omega =\sum\limits_{i=1}^n\frac{A_i(a)}{z-a_i}d(z-a_i) +\sum\limits_{k=1}^n \gamma_k(a)da_k
+\sum\limits_{l=1}^n\sum\limits_{k=1}^n\sum\limits_{m=1}^{r_l}\frac{\gamma_{m, k,  l}(a)}{(z-a_l)^m}da_k,
\end{gather*}
where $\gamma_{m, k, l}(a)$, $\gamma_{k}(a)$ are holomorphic in $D(a^0)$ and $r_l$ is a~maximal~$l$-resonance of system \eqref{syst} for $a=a^0$.
\end{theorem}
We remark that the last terms may be non-zero only if system~\eqref{syst} has some resonant singularities.

The famous example~\cite{ Bol2, Bol7, Bol1} of A.A.~Bolibruch of the non-Schlesinger normalized isomonodromic
deformation is given as follows.

\begin{example}
The family of Fuchsian systems  with four f\/inite singularities
\begin{gather}
\frac{dy}{dz}=\Biggl(\left(
\begin{matrix}
1 & 0
\\
-\frac{2a}{a^2-1} & 0
\end{matrix}
\right) \frac{1}{z+a}+ \left(
\begin{matrix}
0 & -6a
\\
0 & -1
\end{matrix}
\right) \frac{1}{z}
\nonumber
\\
\phantom{\frac{dy}{dz}=}{}
 +\left(
\begin{matrix}
2 & 3+3a
\\
\frac{1}{1+a} & -1
\end{matrix}
\right) \frac{1}{z-1}+ \left(
\begin{matrix}
-3 & 3a-3
\\
\frac{1}{a-1} & 2
\end{matrix}
\right) \frac{1}{z+1} \Biggr)y
\label{ex_nonslez_def}
\end{gather}
is isomonodromic with the dif\/ferential form~$\omega$ given~by
\begin{gather*}
\omega=\left(
\begin{matrix} 1 & 0
\\
-\frac{2a}{a^2-1} & 0
\end{matrix}
\right) \frac{d(z+a)}{z+a}+ \left(
\begin{matrix} 0 & -6a
\\
0 & -1
\end{matrix}
\right) \frac{dz}{z}+ \left(
\begin{matrix} 0 & 0
\\
\frac{2a}{a^2-1} & 0
\end{matrix}
\right) \frac{da}{z+a}
\\
\phantom{\omega=}{}
 +\left(
\begin{matrix} 2 & 3+3a
\\
\frac{1}{1+a} & -1
\end{matrix}
\right) \frac{d(z-1)}{z-1}+ \left(
\begin{matrix} -3 & 3a-3
\\
\frac{1}{a-1} & 2
\end{matrix}
\right) \frac{d(z+1)}{z+1}.
\end{gather*}
The family \eqref{ex_nonslez_def} is isomonodromic by Theorem~\ref{th_isom}.
Here $a_1=-a$ is the resonant singular point (the corresponding residue matrix has eigenvalues $1$ and~$0$, so the
dif\/ference is a~natural number and the maximal 1-resonance for the singular point $a_1=-a$ is equal to~$1$) and $a_2=0$
is also resonant, 2-resonance is equal to~$1$.
The deformation is normalized (since there is no term of the form~$\gamma(a) da$ for any holomorphic $\gamma(a)$) and it
cannot be reduced to any Schlesinger deformation (since there is a~term $\left(
\begin{matrix}
0 & 0
\\
\frac{2a}{a^2-1} & 0
\end{matrix}
\right) \frac{da}{z+a}$ and it cannot be removed by any holomorphic transformation $y=\Gamma(a)y_1$).
\end{example}

We remark that the singularity $z=-a$ in the example above is apparent.
This can be seen by applying a~series of transformations to~\eqref{ex_nonslez_def}.
In particular, changing $y= {\rm diag}(a+z,1)y_1$ gives a~system with a~zero residue matrix at $z=-a$ and, by using an
additional transformation $y_1={\rm diag}(1/z,1)y_2$ we get a~simple system with three singularities and integer
coef\/f\/icients in all residue matrices.

Finally, we remark that in~\cite{Kitaev} A.~Kitaev also considered non-Schlesinger isomonodromic deformations.
He considered some analogues of the Schlesinger equation and studied relations between their solutions and solutions of
some Schlesinger equation.

\section{Middle convolution}

Middle convolution $mc_{\mu}$ is an operation on tuples of residue matrices of a~Fuchsian system introduced~by
M.~Dettweiler and S.~Reiter~\cite{DR,DR2}.
For a~given parameter $\mu \in\mathbb{C}$ one def\/ines residue matrices of dimension $np\times np$ which are partitioned
into blocks and have only one non-zero block consisting of initial residue matrices and the parameter~$\mu$.
By f\/inding invariant subspaces and reducing the size of the matrices (if the invariant subspaces are non-empty) one gets
a~new Fuchsian system with the same singularities but with new residue matrices.
Note that the size of the residue matrices of the resulting system depends on the choice of the parameter~$\mu$.
Middle convolution $mc_{\mu}$ is, in fact, related to the Euler transformation as shown in~\cite{DR2}.

Originally N.~Katz~\cite{Katz} gave an algorithm to construct all physically rigid irreducible local systems on
$\mathbb{CP}^1 \setminus \{\text{f\/inite points}\}$.
M.~Dettweiler and S.~Reiter~\cite{DR,DR2} introduced a~purely algebraic algorithm and showed its equivalence to the
algorithm of Katz.
Although initially the algorithm was applied to rigid systems, as it is shown in~\cite{GF_Painleve, HaraokaFilipuk}, one
can prove non-trivial results for non-rigid systems as well.

Now following~\cite{DR, DR2, GF_Heun, HaraokaFilipuk} we brief\/ly summarize the main steps of the algorithm.
Let ${\bf A}=(A_1,\ldots,A_n)$, $A_k\in\mathbb{C}^{p\times p}$ be the residue matrices of the Fuchsian system
\begin{gather}
\label{mc1}
\frac{dy}{dz}=\sum\limits_{k=1}^n \frac{A_k}{z-a_k}y
\end{gather}
with f\/inite singular points $a_k\in\mathbb{C}$, $k=1,\ldots,n$.

For a~given parameter $\mu\in\mathbb{C}$ one def\/ines the so-called convolution matrices ${\bf B}=c_{\mu}({\bf
A})=(B_1,\ldots,B_n)$~by
\begin{gather*}
B_k=\left(
\begin{matrix}
0 & \dots & 0 & 0 & 0 & \dots & 0
\\
\vdots & &\vdots &\vdots &\vdots & &\vdots
\\
A_1 & \dots & A_{k-1} & A_k+\mu I_p & A_{k+1} & \dots & A_n
\\
\vdots & &\vdots &\vdots &\vdots & &\vdots
\\
0 & \dots & 0 & 0 & 0 & \dots & 0
\end{matrix}
\right)\in\mathbb{C}^{np\times np}
\end{gather*}
such that $B_k$ is zero outside the~$k$-th block row.

By using the convolution matrices we def\/ine a~new Fuchsian system of order $np$ with the same number of singularities as
in the original system:
\begin{gather*}
\frac{dy_1}{dz}=\sum\limits_{k=1}^n \frac{B_k}{z-a_k}y_1.
\end{gather*}
There are the following invariant subspaces (with respect to $B_k$) of the column vector space $\mathbb{C}^{np}$:
\begin{gather*}
\mathcal{L}_k=\left(
\begin{matrix}
0
\\
\vdots
\\
0
\\
\rm{Ker}(A_k)
\\
0
\\
\vdots
\\
0
\end{matrix}
\right) (k {\text{-th~entry})},\qquad k=1,\ldots,n,
\end{gather*}
and
\begin{gather*}
\mathcal{K}=\bigcap_{k=1}^{n}{\rm{Ker}}(B_k)={\rm{Ker}}(B_1+\cdots+B_n).
\end{gather*}
Next we f\/ix an isomorphism between $\mathbb{C}^{np}/(\mathcal{K+L})$, where $\mathcal{L}=\oplus_{k=1}^n\mathcal{L}_k$,
and $\mathbb{C}^m$ for some~$m$.
Note that if there are non-zero invariant subspaces, then $m<np$.
\begin{definition}
The matrices $\widetilde{{\bf{B}}}=mc_\mu({\bf{A}}):=(\widetilde{B}_1,\ldots,\widetilde{B}_n)\in\mathbb{C}^{m\times m}$,
where $\widetilde{B}_k$ is induced by the action of $B_k$ on $\mathbb{C}^m \simeq \mathbb{C}^{np}/(\mathcal{K+ L})$, are
called {\it the additive version of the middle convolution of~${\bf A}$ with the parameter~$\mu$}.
\end{definition}

Finally, the resulting Fuchsian system of order~$m$ is given~by
\begin{gather*}
\frac{dy_2}{dz}=\sum\limits_{k=1}^n \frac{\widetilde{B}_k}{z-a_k}y_2.
\end{gather*}

Thus, middle convolution $mc_{\mu}$ is a~transformation on tuples of matrices
\begin{gather*}
((A_1,\ldots,A_n)\in\bigl(\mathbb{C}^{p\times p}\bigr)^n)\rightarrow
mc_{\mu}(A_1,\ldots,A_n)=\big(\widetilde{B}_1,\ldots,\widetilde{B}_n\big)\in \bigl(\mathbb{C}^{m\times m}\bigr)^n.
\end{gather*}

In~\cite{HaraokaFilipuk} the authors were interested whether middle convolution preserves Schlesinger's deformation
equations~\eqref{Schl eqs} for non-resonant Fuchsian systems and gave an af\/f\/irmative answer to this question.
In particular, the following statement holds.

\begin{theorem}[\cite{HaraokaFilipuk}]
\label{th_HF}
Assume that there is no integer differences between any two distinct eigenvalues of $A_j$ and the Jordan canonical form
of $A_j$ is independent of $a_1, a_2,\ldots, a_n$ for $j = 1, 2,\ldots, n + 1$.
Then the Hitchin systems for the Fuchsian systems obtained by addition and middle convolution with parameters
independent of $a_1, a_2,\ldots, a_n$ coincide with the Hitchin system for the initial Fuchsian system~\eqref{mc1}.
\end{theorem}

Note that Theorem~\ref{th_HF} cannot be applied to a~system with a~resonant Fuchsian singular point, so it cannot be
applied to non-Schlesinger isomonodromic deformations.

\section{Non-Schlesinger isomonodromic deformations\\ and middle convolution for resonant Fuchsian systems}

Since middle convolution preserves Schlesinger isomonodromic deformation equations~\cite{HaraokaFilipuk}, it is natural
to ask what happens to non-Schlesinger isomonodromic deformations under middle convolution.
In this paper we construct an explicit example to show that middle convolution does not in general preserve
non-Schlesinger isomonodromic deformations.
Our examples are based on the modif\/ications of the Bolibruch example.
Note that it is very important and of interest to specialists in the deformation theory to f\/ind new explicit examples of
non-Schlesinger iso\-monodro\-mic deformations because of the dif\/f\/iculty to write down the corresponding dif\/ferential
1-form~$\omega$~\cite{Bibilo_Thesis}.
Therefore, any new examples especially for Fuchsian systems of order higher than 2 are necessary.
Our explicit examples show that under middle convolution the resonance condition may appear or disappear.
Moreover, the maximal~$i$-resonance of the system may change.
This shows that in some cases the non-Schlesinger isomonodromic deformations are not preserved by middle convolution.

{\bf From non-Schlesinger to non-Schlesinger deformations}.
We shall present an example illustrating how to obtain new non-Schlesinger isomonodromic deformations by using middle
convolution with changing the order of the resulting Fuchsian system.
We remark that this example is non-trivial and new.
We underline that the main dif\/f\/iculty with explicit examples of the non-Schlesinger deformations is to f\/ind the relevant
dif\/ferential 1-form.

In this illustrative example we prove that after applying middle convolution with a~parameter~$\mu$, $\mu\neq 0$, to
a~certain non-Schlesinger $(2\times2)$-isomonodromic family we get a~new $(5\times 5)$-family which is non-Schlesinger
and it cannot be reduced to any Schlesinger isomonodromic deformation.

Consider a~family of Fuchsian systems
\begin{gather}
\frac{dy}{dz}=\Biggl(\left(
\begin{matrix} \frac{1}{2} & 0
\\
-\frac{2a}{a^2-1} & -\frac{1}{2}
\end{matrix}
\right) \frac{1}{z+a}+ \left(
\begin{matrix} 0 & -6a
\\
0 & -\frac{1}{2}
\end{matrix}
\right) \frac{1}{z}
\nonumber
\\
\phantom{\frac{dy}{dz}=}{}
 +\left(
\begin{matrix} 1 & 3a+3
\\
\frac{1}{1+a} & 3
\end{matrix}
\right) \frac{1}{z-1}+ \left(
\begin{matrix} -\frac{3}{2} & 3a-3
\\
\frac{1}{a-1} & -2
\end{matrix}
\right) \frac{1}{z+1} \Biggr)y.
\label{ex2}
\end{gather}
This family has four singular points (since the residue matrix at inf\/inity vanishes) with two resonant singularities:
$a_1=-a$ with the maximal 1-resonance equal to $1$ and $a_3=1$ with the maximal 3-resonance equal to~$4$.
The family is isomonodromic since the conditions of Theorem~\ref{th_isom} are fulf\/illed with the dif\/ferential form given~by
\begin{gather*}
\omega= \left(
\begin{matrix} \frac{1}{2} & 0
\\
-\frac{2a}{a^2-1} & -\frac{1}{2}
\end{matrix}
\right) \frac{d(z+a)}{z+a}+ \left(
\begin{matrix} 0 & -6a
\\
0 & -\frac{1}{2}
\end{matrix}
\right) \frac{dz}{z}+ \left(
\begin{matrix} 1 & 3a+3
\\
\frac{1}{1+a} & 3
\end{matrix}
\right) \frac{d(z-1)}{z-1}
\\
\phantom{\omega=}{}
 +\left(
\begin{matrix} -3/2 & 3a-3
\\
\frac{1}{a-1} & -2
\end{matrix}
\right) \frac{d(z+1)}{z+1} + \left(
\begin{matrix} 0 & 0
\\
\frac{2a}{a^2-1} & 0
\end{matrix}
\right) \frac{da}{z+a}.
\end{gather*}

We remark that similarly to Bolibruch's example, the family~\eqref{ex2} has an apparent singularity at $z=-a$ and can be
transformed to a~system with four singular points, including the point at inf\/inity, and with residue matrices which do
not depend on the parameter.

By applying middle convolution with the parameter~$\mu$ ($\mu \notin \frac{1}{2} \mathbb{Z}$) to this family of Fuchsian
systems and adding the vectors $e_1$, $e_2$, $e_4$, $e_6$, $e_8$ to the basis of the invariant subspaces we get a~new $(5\times
5)$-family.
It has one resonant Fuchsian singularity $a_1=-a$ with the maximal 1-resonance equal to 1.
We have
\begin{gather*}
\frac{dy}{dz}=\left(\frac{A_1(a)}{z+a}+\frac{A_2(a)}{z}+\frac{A_3(a)}{z-1}+\frac{A_4(a)}{z+1} \right) y,
\end{gather*}
where
\begin{gather*}
 A_1(a)=\left(
\begin{matrix}
\mu+1/2 & 0& -6a& 3a+3& 3a-3
\\
-\frac{2 a}{a^2-1}& \mu-\frac{1}{2}& -\frac{1}{2}& 3& -2
\\
0& 0& 0& 0& 0
\\
0& 0& 0& 0& 0
\\
0& 0& 0& 0& 0
\end{matrix}
\right),
\\
 A_2(a)=\left(
\begin{matrix}
0& 0& 0& 0& 0
\\
0& 0& 0& 0& 0
\\
-\frac{2a}{-1+a^2}& -\frac{1}{2}& \mu-\frac{1}{2}& 3& -2
\\
0& 0& 0& 0& 0
\\
0& 0& 0& 0& 0
\end{matrix}
\right),
\\
 A_3(a)=\left(
\begin{matrix}
0& 0& 0& 0& 0
\\
0& 0& 0& 0& 0
\\
0& 0& 0& 0& 0
\\
-\frac{11a+1}{6(a^2-1)}&-\frac{1}{2}& -\frac{5a+1}{2(a+1)}& \mu+4& -\frac{a+3}{a+1}
\\
0& 0& 0& 0& 0
\end{matrix}
\right),
\\
 A_4(a)=\left(
\begin{matrix}
0& 0& 0& 0& 0
\\
0& 0& 0& 0& 0
\\
0& 0& 0& 0& 0
\\
0& 0& 0& 0& 0
\\
-\frac{9a+1}{4(a^2-1)}& -\frac{1}{2}& \frac{5a+1}{2(a-1)}&\frac{3(a-3)}{2(a-1)}& \mu-\frac{7}{2}
\end{matrix}
\right).
\end{gather*}

This Fuchsian family with f\/ive singular points is isomonodromic since it is def\/ined by the dif\/ferential 1-form which
satisf\/ies conditions of Theorem~\ref{th_isom}.
This form is given~by
\begin{gather*}
\omega=A_1(a)\frac{d(z+a)}{z+a}+A_2(a)\frac{dz}{z}+A_3(a)\frac{d(z-1)}{z-1}+A_4(a)\frac{d(z+1)}{z+1}
\\
\phantom{\omega=}{}
 +\left(
\begin{matrix}
0& 0& 6& -3& -3
\\
\frac{2(2a^2+5a+1)}{(a^2-1)^2}& \frac{2a^2-5a+1}{2(a^2-1) a}& -\frac{47 a^2+1}{2(a^2-1)a}& \frac{9a+3}{a^2-1}&
\frac{2(7a+1)}{a^2-1}
\\
0& \frac{1}{2a}& -\frac{1+5a}{2a-2a^3}& 0& 0
\\
-\frac{1}{6(1+a)^2}& \frac{1}{2(a+1)}& 0& \frac{3}{a^2-1}& 0
\\
\frac{1}{4(a-1)^2}& \frac{1}{2(a-1)}& 0& 0& \frac{2}{a^2-1}
\end{matrix}
\right)da
\\
\phantom{\omega=}{}
 +\left(
\begin{matrix}
0& 0& 0& 0& 0
\\
-\frac{2(2\mu+1)a}{a^2-1}& 0& \frac{24a^2}{a^2-1}& -\frac{12a}{a-1}& -\frac{12a}{a-1}
\\
0& 0& 0& 0& 0
\\
0& 0& 0& 0& 0
\\
0& 0& 0& 0& 0
\end{matrix}
\right) \frac{da}{z+a}.
\end{gather*}
Moreover, the last term of the form $\psi(a)\frac{da}{z+a}$ shows that this isomonodromic deformation is non-Schlesinger
and cannot be reduced to any Schlesinger deformation.
Also it is non-normalized because of the non-zero term of the form $\gamma(a)da$.

We also remark that in general it is computationally dif\/f\/icult to f\/ind the non-Schlesinger isomonodromic dif\/ferential
1-form for a~given resonant Fuchsian family such that Theorem~\ref{th_isom} is fulf\/illed.
Actually, the form in our example is not unique, there exists a~one-parameter family of such forms but we presented only
the case when the parameter is zero.
It is possible to simplify the calculations a~bit by f\/irst deriving some identities from $d\omega=\omega\wedge\omega$,
where~$\omega$ is a~1-form with holomorphic coef\/f\/icients, i.e., $\omega=A_1(a) dz/(z+a)+A_2(a)
dz/z+A_3(a)dz/(z-1)+A_4(a)dz/(z+1)+A_5(a) da/(z+a)+B(a) da$.
For example, we have $[A_1,A_5]=A_1-A_5$, $\partial(A_{\infty})/\partial a=-[A_{\infty},B]$, where
$A_{\infty}=-(A_1+A_2+A_3+A_4)$, and so on.

{\bf Example showing that non-Schlesinger deformations are not preserved by middle convolution.} In this example we
apply middle convolution with a~parameter~$\mu$, $\mu \neq 0$, to a~certain non-Schlesinger isomonodromic $(2\times
2)$-family (which is a~modif\/ied Bolibruch examp\-le~\eqref{ex_nonslez_def}) with f\/ive singular points (two of them are
resonant).
In the resulting family the resonances disappear and, hence, the resulting Fuchsian system cannot be included in any
Bolibruch's non-Schlesinger isomonodromic deformation with a~non-zero third term by Theo\-rems~\ref{th_isom} and~\ref{th 2.2}.

Consider a~family
\begin{gather}
\frac{dy}{dz}=\Biggl(\left(
\begin{matrix} 1 & 0
\\
-\frac{2a}{a^2-1} & 0
\end{matrix}
\right) \frac{1}{z+a}+ \left(
\begin{matrix} 0 & -6a
\\
0 & -1
\end{matrix}
\right) \frac{1}{z}
\nonumber
\\
\phantom{\frac{dy}{dz}=}{}
+\left(
\begin{matrix} \frac{3}{2}-\frac{\sqrt{21}}{2} & 3+3a
\\
\frac{1}{1+a} & -\frac{3}{2}-\frac{\sqrt{21}}{2}
\end{matrix}
\right) \frac{1}{z-1}+ \left(
\begin{matrix} -3 & 3a-3
\\
\frac{1}{a-1} & 2
\end{matrix}
\right) \frac{1}{z+1} \Biggr)y
\label{ex3}
\end{gather}
with f\/ive singular points $-a$, $0$, $1$, $-1$, $\infty$.
The system has two resonant Fuchsian points $a_1=-a$, $a_2=0$ with the maximal resonances both equal to~$1$.
Note that this family is obtained by addition from Bolibruch's example.
Recall that an addition for a~Fuchsian system is a~transformation of the form $y\rightarrow (z-a_i)^{-\alpha}y$ such
that the corresponding to the singularity $a_i$ residue matrix $A_i$ is transformed by a~shift, i.e., $A_i\rightarrow
A_i+\alpha I_p$, where $I_p$ is an identity matrix of the same dimension as $A_i$.
The residue matrix at inf\/inity is the scalar.
Therefore, it is straightforward to construct a~corresponding 1-form.
In particular, this family is isomonodromic since the conditions in Theorem~\ref{th_isom} are fulf\/illed with the
dif\/ferential form given~by
\begin{gather*}
\omega=\left(
\begin{matrix} 1 & 0
\\
-\frac{2a}{a^2-1} & 0
\end{matrix}
\right) \frac{d(z+a)}{z+a}+ \left(
\begin{matrix} 0 & -6a
\\
0 & -1
\end{matrix}
\right) \frac{dz}{z}+ \left(
\begin{matrix} \frac{3}{2}-\frac{\sqrt{21}}{2} & 3+3a
\\
\frac{1}{1+a} & -\frac{3}{2}-\frac{\sqrt{21}}{2}
\end{matrix}
\right) \frac{d(z-1)}{z-1}
\\
\phantom{\omega=}{}
 +\left(
\begin{matrix} -3 & 3a-3
\\
\frac{1}{a-1} & 2
\end{matrix}
\right) \frac{d(z+1)}{z+1}+ \left(
\begin{matrix}
0 & 0
\\
\frac{2a}{a^2-1} & 0
\end{matrix}
\right) \frac{da}{z+a}.
\end{gather*}

Applying middle convolution with the parameter $\mu=(1+\sqrt{21})/2$, which is chosen to be equal to the eigenvalue of
the residue matrix at inf\/inity, and adding vectors $e_7$, $e_1-e_4$ to the basis of the invariant subspaces, we get a~new
$(3\times3)$-family.
The resulting family has no resonant singular points, so it cannot be an isomonodromic deformation with the dif\/ferential
form having the non-zero term of the form $\psi(a) \frac{da}{z+a}$.
One can also calculate that it is not a~Schlesinger isomonodromic deformation.
We were not able to f\/ind a~f\/lat dif\/ferential 1-form
$\omega=\frac{A_1(a)d(z+a)}{z+a}+\frac{A_2(a)dz}{z}+\frac{A_3(a)dz}{z-1}+\frac{A_4(a)dz}{z+1}+X(a)da$ for
a~holomorphic~$X$ to show that it is a~non-Schlesinger isomonodromic deformation, which can be transformed to the Schlesinger one.
Thus, middle convolution may destroy resonances in the residue matrices of the Fuchsian systems and, hence,
non-Schlesinger deformations, which cannot be reduced to the Schlesinger ones.
Explicitly, the resulting family of the Fuchsian systems is
\begin{gather}
\label{ex3 mc}
\frac{dy}{dz}=\left(\frac{A_1}{z+a}+\frac{A_2}{z}+\frac{A_3}{z-1}+\frac{A_4}{z+1}\right)y,
\end{gather}
where
\begin{gather*}
A_1=\left(
\begin{matrix} -6 a~& -3 & \frac{1}{2}(12a+\sqrt{21}+3)
\\
0& 0& 0
\\
-6 a& -3& \frac{1}{2} (12 a+ \sqrt{21}+3)
\end{matrix}
\right),
\\
A_2=\left(
\begin{matrix} \frac{1}{2} (\sqrt{21}-1) & \frac{1}{a-1}& - \frac{2a}{a^2-1} - \frac{\sqrt{21}}2+\frac{1}{2}
\\
0& 0& 0
\\
0& 0& 0
\end{matrix}
\right),
\\
A_3=\left(
\begin{matrix}
\frac{1}{2} ((9 - \sqrt{21})a - \sqrt{21}-3)& \frac{(\sqrt{21}+3)(a+1)}{2 (a-1)}+3& \alpha
\\
\frac{1}{2} ((9 - \sqrt{21})a- \sqrt{21}-3)& \frac{(\sqrt{21}+3)(a+1)}{2 (a-1)}+3& \alpha
\\
\frac{1}{2} ((9 - \sqrt{21})a - \sqrt{21}-3)& \frac{(\sqrt{21}+3)(a+1)}{2 (a-1)}+3& \alpha
\end{matrix}
\right),
\\
\alpha=\frac{1}{2} (\sqrt{21}+3) (a+1)-6 a-1 -\frac{(\sqrt{21}+3)a}{a-1},
\\
A_4=\left(
\begin{matrix}
\frac{1}{2}((\sqrt{21}+3)a+\sqrt{21}+5) & -\frac{(\sqrt{21}+3) (a+1)+2}{2(a-1)} & -
\frac{((\sqrt{21}+3)a+\sqrt{21}+5)(a^2-2a-1)}{2(a^2-1)}
\\
\frac{1}{2}((\sqrt{21}-9)a+\sqrt{21}+3)& -\frac{4a+\sqrt{21}-1}{(a-1)}& -\frac{(\sqrt{21}-9)
a^2-2(\sqrt{21}-2)a-\sqrt{21}-1}{2(a-1)}
\\
\frac{1}{2} (\sqrt{21}+3) (a+1)& - \frac{(\sqrt{21}+3) (a+1)}{2 (a-1)}& - \frac{(\sqrt{21}+3) (a^2-2a-1)}{2 (a-1)}
\end{matrix}
\right).
\end{gather*}

It is interesting to note that family~\eqref{ex3 mc} can be transformed to a~resonant one by a~series of special
transformations.
The eigenvalues of the matrix $A_1$ are $\mu+1$, $0$, $0$, where $\mu $ is the parameter of middle convolution def\/ined above,
for $A_2$ the eigenvalues are $\mu-1$, $0$, $0$ and the matrix at inf\/inity is diagonal with all diagonal elements equal to~$\mu$.
If in family~\eqref{ex3 mc} we f\/irst shift the residue matrix $A_1 \rightarrow A_1-\mu I_3$ and then in the resulting
family replace
\begin{gather*}
%\label{tr}
y=\left(
\begin{matrix}
(z+a)&f_1(a)&f_2(a)
\\
0&1&\frac{f_2(a)-1}{f_1(a)}
\\
0&0&1
\end{matrix}
\right)y_1,
\end{gather*}
where $f_1(a)$ and $f_2(a)$ are arbitrary holomorphic functions, we get a~new family for $y_1$ with a~resonant Fuchsian
singularity at $z=-a$ and a~resonant Fuchsian singularity at inf\/inity.
Moreover, applying the transformation $y_1={\rm diag}(1/z,1,1)y_2$, we can get rid of the singularity at inf\/inity.
The resulting Fuchsian system remains resonant.
See also~\cite{Pober} for other special transformations of Fuchsian systems.

\section{Discussion and open problems}

In this paper we have shown that middle convolution is useful in constructing new explicit examples of non-Schlesinger
isomonodromic deformations for resonant Fuchsian systems.
We have also shown that after middle convolution non-Schlesinger isomonodromic deformations may stay (see example~\eqref{ex2}) or not (see example~\eqref{ex3}).
It is a~new fact about middle convolution.

There are many open questions and problems in the theory of resonant Fuchsian systems and in understanding and
explaining the results of middle convolution for them.
They include:
\begin{enumerate}\itemsep=0pt
\item The description of suf\/f\/icient conditions when non-Schlesinger isomonodromic deformations of a~Fuchsian system with
a~resonant singularity stay after middle convolution.
\item The description of suf\/f\/icient conditions of the existence of non-Schlesinger isomonodromic deformations which
cannot be reduced to the Schlesinger ones.
\item What kind of isomonodromic dif\/ferential 1-form is obtained after middle convolution in general? How to reduce the
complexity of straightforward calculations of the iso\-mo\-nodro\-mic dif\/ferential 1-form?
\item Is there any analogue of the Hitchin systems for resonant Fuchsian systems which will be invariant under middle
convolution?
\item The construction of more examples of non-Schlesinger isomonodromic deformations (with or without the application
of middle convolution, with and without apparent singularities).
In particular, to our knowledge, there are no explicit examples of non-Schlesinger iso\-monodro\-mic deformations with the
order of poles greater then 1 in the ``non-Schlesinger'' term of the isomonodromic form.
\item The (practical) description how the maximal resonance changes after middle convolution, including the cases when
the resulting system is non-resonant.
\end{enumerate}

\looseness=-1
Non-Schlesinger isomonodromic deformations have a~natural analogue in theory of meromorphic systems with irregular
singularities~\cite{Bibilo}.
As mentioned in the introduction, there also exist generalizations of middle convolution to irregular systems
(e.g.,~\cite{Hiroe, Kawakami, Takemura_irreg}).
Therefore, almost all questions and problems mentioned above can be formulated for irregular systems as well.
In particular, it is desirable to understand middle convolution for both resonant and non-resonant systems with
irregular singularities and, in particular, to f\/ind more explicit examples of non-Schlesinger isomonodromic deformations
of meromorphic systems with irregular resonant singularities.

\subsection*{Acknowledgments}

Part of this work was carried out while Yu.~Bibilo was visiting the Univerity of Warsaw in April 2014.
The authors acknowledge the support of the Polish NCN Grant 2011/03/B/ST1/00330.
Yu.~Bibilo also acknowledges the support of the Russian Foundation for Basic Research (grant no.\
RFBR 14-01-00346~A).

\pdfbookmark[1]{References}{ref}
\LastPageEnding

\end{document}